\documentclass[10pt]{article}
\font\smallit=cmti10
\font\smalltt=cmtt10
 
\usepackage{amssymb,latexsym,amsmath,epsfig,amsthm}

\makeatletter

\renewcommand\section{\@startsection {section}{1}{\z@}
{-30pt \@plus -1ex \@minus -.2ex}
{2.3ex \@plus.2ex}
{\normalfont\normalsize\bfseries\boldmath}}

\renewcommand\subsection{\@startsection{subsection}{2}{\z@}
{-3.25ex\@plus -1ex \@minus -.2ex}
{1.5ex \@plus .2ex}
{\normalfont\normalsize\bfseries\boldmath}}

\renewcommand{\@seccntformat}[1]{\csname the#1\endcsname. }

\makeatother
\theoremstyle{plain}

\newtheorem{proposition}{Proposition}
\newtheorem{corollary}{Corollary}
\newtheorem{lemma}{Lemma}
\newtheorem{definition}{Definition}

\theoremstyle{remark}

\def\la{\langle}
\def\ra{\rangle}
\def\be{\begin{equation}}
\def\ee{\end{equation}}
\def\bea{\begin{eqnarray}}
\def\eea{\end{eqnarray}}

\begin{document}
\begin{center}
\uppercase{\bf Symmetric (not Complete Intersection) Numerical Semigroups 
Generated by Five Elements}
\vskip 20pt
{\bf Leonid G. Fel\footnote{The research was partly supported by the Kamea 
Fellowship.}}\\
{\smallit Department of Civil Engineering, Technion, Haifa 32000, Israel}\\
{\tt lfel@technion.ac.il}\\
\end{center}
\vskip 10pt
\centerline{\smallit Received: , Revised: , Accepted: , Published: }
\vskip 30pt 
\centerline{\bf Abstract}

\noindent
We consider symmetric (not complete intersection) numerical semigroups
${\sf S}_5$, generated by five elements, and derive inequalities for degrees
of syzygies of ${\sf S}_5$ and find the lower bound $F_5$ for their Frobenius
numbers. We study a special case ${\sf W}_5$ of such semigroups, which satisfy
the Watanabe Lemma, and show that the lower bound $F_{5w}$ for the Frobenius 
number of the semigroup ${\sf W}_5$ is stronger than $F_5$.
\pagestyle{myheadings}
\markright{\smalltt INTEGERS: 18 (2018)\hfill}
\thispagestyle{empty}
\baselineskip=12.875pt
\vskip 30pt 
\section{Symmetric numerical semigroups generated by five integers}\label{s1}
Let a numerical semigroup ${\sf S}_m=\la d_1,\ldots,d_m\ra$ be generated by a
set of $m$ integers $\{d_1,\ldots,d_m\}$ such that $\gcd(d_1,\ldots,d_m)=1$,
where $d_1$ and $m$ denote a multiplicity and an embedding dimension ({\em 
edim}) of ${\sf S}_m$. There exist $m-1$ polynomial identities \cite{fel017} 
for degrees of syzygies associated with the semigroup ring $k[{\sf S}_m]$. They 
are a source of various relations for semigroups of a different nature. In the 
case of complete intersection (CI) semigroups, such a relation for the degrees 
$e_j$ of the 1st syzygy was found in \cite{fel017}. The next nontrivial case 
exhibits a symmetric (not CI) semigroup generated by $m\ge 4$ integers. In 
\cite{fel015} such semigroups with $m=4$ were studied and the lower bound for 
the Frobenius number $F({\sf S}_4)$ was found. In the present paper we deal 
with a more difficult case of symmetric semigroups ${\sf S}_5$.

Consider a symmetric numerical semigroup ${\sf S}_5$, which is not CI, and
generated by five elements $d_i$. Its Hilbert series $H\left({\sf S}_5;t
\right)$ with the 1st Betti number $\beta_1$ reads,
\bea
&&H\left({\sf S}_5;t\right)=\frac{Q_5(t)}{\prod_{i=1}^5\left(1-t^{d_i}\right)},
\nonumber\\
&&Q_5(t)=1-\sum_{j=1}^{\beta_1}t^{x_j}+\sum_{j=1}^{\beta_1-1}\left(t^{y_j}+
t^{g-y_j}\right)-\sum_{j=1}^{\beta_1}t^{g-x_j}+t^{g},\label{h1}
\eea
where $x_j,y_j,g\in{\mathbb Z}_{>0}$, $x_j,y_j<g$, and the Frobenius number is
defined as follows: $F({\sf S}_5)=g-\sigma_1$, where $\sigma_1=\sum_{j=1}^5d_j$.
There are two more constraints, $\beta_1>4$ and $d_1>5$. The inequality $\beta_1
>4$ holds since the semigroup ${\sf S}_5$ is not CI, and the condition $d_1>5$ 
is necessary since the numerical semigroup $\la m,d_2,\ldots,d_m\ra$ is never 
symmetric \cite{fel11}.

Polynomial identities for degrees of syzygies for numerical semigroups were
derived in \cite{fel017}, Theorem 1. In the case of the symmetric (not CI) 
semigroup ${\sf S}_5$, they read:
\bea
&&\sum_{j=1}^{\beta_1}x_j^r-\sum_{j=1}^{\beta_1-1}\left[y_j^r+(g-y_j)^r
\right]+\sum_{j=1}^{\beta_1}(g-x_j)^r-g^r=0,\quad 1\le r\le 3,\label{h2}\\
&&\sum_{j=1}^{\beta_1}x_j^4-\sum_{j=1}^{\beta_1-1}\left[y_j^4+(g-y_j)^4\right]
+\sum_{j=1}^{\beta_1}(g-x_j)^4-g^4=-24\pi_5,\quad\pi_5=\prod_j^5d_j.\qquad
\label{h3}
\eea
Only two of the four identities in (\ref{h2},\ref{h3}) are not trivial; these 
are the identities with $r=2,4$,
\bea
&&\sum_{j=1}^{\beta_1}x_j(g-x_j)=\sum_{j=1}^{\beta_1-1}y_j(g-y_j),\label{h4}\\
&&\sum_{j=1}^{\beta_1}x_j^2(g-x_j)^2+12\pi_5=\sum_{j=1}^{\beta_1-1}y_j^2(g-
y_j)^2.\label{h5}
\eea
Denote by $Z_k$ the k-th power symmetric polynomial $Z_k(z_1,\ldots,z_n)=\sum_{
j=1}^nz_j^k$, $z_j\ge 0$, and recall the Newton-Maclaurin inequality 
\cite{har59},
\bea
Z_1^2\le nZ_2.\label{h6}
\eea
Applying (\ref{h6}) to the right-hand side of equality (\ref{h5}) and 
substituting into the resulting inequality an equality (\ref{h4}), we obtain
\be
\left(\sum_{j=1}^{\beta_1}x_j(g-x_j)\right)^2\le (\beta_1-1)\left[\sum_{j=1}^
{\beta_1}x_j^2(g-x_j)^2+12\pi_5\right].\label{h7}
\ee
On the other hand, according to (\ref{h6}), another inequality holds:
\bea
\left(\sum_{j=1}^{\beta_1}x_j(g-x_j)\right)^2\le\beta_1\sum_{j=1}^{\beta_1}x_j
^2(g-x_j)^2.\label{h8}
\eea
\begin{proposition}\label{pr1}
Let a symmetric (not CI) semigroup ${\sf S}_5$ be given with the Hilbert series
$H\left({\sf S}_5;z\right)$ according to (\ref{h1}). Then the following
inequality holds:
\be
g\ge g_5,\qquad g_5=\lambda(\beta_1)\sqrt[4]{\pi_5},\qquad\lambda(\beta_1)=
4\sqrt[4]{\frac{3}{4}\frac{\beta_1-1}{\beta_1}}.\label{h9}
\ee
\end{proposition}
\begin{proof}
Inequality (\ref{h8}) holds always, while inequality (\ref{h7}) is not valid for
every set $\{x_1,\ldots,x_{\beta_1},g\}$. In order to make both inequalities
consistent, we have to find a range for $g$ where both inequalities (\ref{h7}) 
and (\ref{h8}) are satisfied. As we will see in (\ref{h12}), this requirement 
also constrains the admissible values of $x_j$. To provide both inequalities to 
be correct, it is enough to require
\be
\frac1{\beta_1-1}\left(\sum_{j=1}^{\beta_1}x_j(g-x_j)\right)^2-12\pi_5\ge
\frac1{\beta_1}\left(\sum_{j=1}^{\beta_1}x_j(g-x_j)\right)^2.\label{h10}
\ee
Making use of the notation in (\ref{h6}), rewrite inequality (\ref{h10}) for 
$X_1=\sum_{j=1}^{\beta_1}x_j$ and $X_2=\sum_{j=1}^{\beta_1}x_j^2$,
\bea
gX_1-C\ge X_2,\qquad C=\sqrt{12\beta_1(\beta_1-1)\pi_5},\nonumber
\eea
and combine it with inequality (\ref{h6}), namely, $X_2\ge \beta_1^{-1}X_1^2$. 
Thus, we obtain
\bea
X_1^2-\beta_1gX_1+\beta_1C\le 0.\label{h11}
\eea
Represent the last inequality (\ref{h11}) in the following way,
\bea
\left(\frac{\beta_1g}{2}\right)^2\ge\beta_1C+\left(X_1-\frac{\beta_1g}{2}
\right)^2,\nonumber
\eea
and obtain immediately the lower bound $g_5$ in accordance with (\ref{h9}).
\end{proof}
The lower bound for the Frobenius number is given by $F_5=g_5-\sigma_1$. 
Inequality (\ref{h11}) constrains the degrees $x_j$ of the 1st syzygy for the 
symmetric (not CI) semigroup ${\sf S}_5$,
\bea
\frac{\beta_1g}{2}\left(1-\sqrt{1-\frac{g_5^2}{g^2}}\right)\le X_1\le
\frac{\beta_1g}{2}\left(1+\sqrt{1-\frac{g_5^2}{g^2}}\right).\label{h12}
\eea
Below we present twelve symmetric (not CI) semigroups ${\sf S}_5$ with different
Betti's numbers $\beta_1=7,8,9,13$ and give a comparative Table 1 for their 
largest degree $g$ of syzygies and its lower bound $g_5$:

\bea
&&\beta_1=7,\quad A_1^7=\la 6,10,14,15,19\ra,\quad A_2^7=
\la 6,10,14,17,21\ra,\quad A_3^7=\la 9,10,11,13,17\ra,\quad\nonumber\\
&&\beta_1=8,\quad A_1^8=\la 6,10,14,19,23\ra,\quad A_2^8=
\la 8,10,13,14,19\ra,\quad A_3^8=\la 8,9,12,13,19\ra,\nonumber\\
&&\beta_1=9,\quad A_1^9=\la 7,12,13,18,23\ra,\quad A_2^9=
\la 9,12,13,14,19\ra,\quad A_3^9=\la 8,11,12,15,25\ra,\nonumber\\
&&\beta_1=13,\;A_1^{13}=\la 19,23,29,31,37\ra,\;A_2^{13}=\la 19,27,28,31,32\ra,
\;A_3^{13}=\la 23,28,32,45,54\ra.\nonumber
\eea

$$
\begin{array}{|c||c|c|c||c|c|c||c|c|c||c|c|c|c|c|c|c|}\hline
{\sf S}_5 & A_1^7 & A_2^7 & A_3^7 & A_1^8 & A_2^8 & A_3^8 & A_1^9 & A_2^9 & 
A_3^9 & A_1^{13} & A_2^{13} & A_3^{13} \\ \hline\hline
g & 87 & 93 & 85 & 99 & 89 & 84 & 102 & 96 & 100 & 240 & 236 & 331 \\ \hline
g_5 & 79.2 & 83.8 & 77.5 & 88.6 & 82.6 & 77.4 & 93.7 & 89.4 & 90.7 & 225.3 &
224.2 & 306.9 \\ \hline
\end{array}
$$

\begin{center}
Table 1. The largest degree $g$ of syzygies and its lower bound $g_5$ for 
symmetric (not CI) semigroups ${\sf S}_5$ with different Betti's numbers 
$\beta_1=7,8,9,13$.
\end{center}

The semigroups $A_i^{13}$, $i=1,2,3$, were studied by H. Bresinsky \cite{bre79};
the other semigroups were generated numerically with the help of the package
"NumericalSgps" \cite{dgm00} by M. Delgado.

Comparing $F_5$ with the two known lower bounds of the Frobenius numbers $F_{CI_
5}$ and $F_{NS_5}$ for the symmetric CI \cite{fel017} and nonsymmetric 
\cite{kil00} semigroups generated by five elements, respectively, we have
\be
F_{CI_5}=4\sqrt[4]{\pi_5}-\sigma_1,\qquad F_{NS_5}=\sqrt[4]{24\pi_5}-\sigma_1,
\qquad F_{NS_5}<F_5<F_{CI_5}.\nonumber
\ee
\section{Symmetric numerical (not CI) semigroups ${\sf S}_5$ with $W$ property}
\label{s2}
Watanabe \cite{wata73} gave a construction of the numerical semigroup $S_m$
generated by $m$ elements starting with a semigroup $S_{m-1}$ generated by 
$m-1$ elements and proved the following lemma.
\begin{lemma}\label{le1}{\rm (\cite{wata73}).}
Let a numerical semigroup $S_{m-1}\!=\!\la\delta_1,\ldots,\delta_{m-1}\ra$ be 
given and $a\in{\mathbb Z}_{>0}$, $a>1$, $d_m>m$, such that $\gcd(a,d_m)=1$, 
$d_m\in S_{m-1}$. Consider a numerical semigroup $S_m\!=\!\la a\delta_1,\ldots,
a\delta_{m-1},d_m\ra$ and denote it by $S_m\!=\!\la aS_{m-1},d_m\ra$. Then $S_m$
is symmetric if and only if $S_{m-1}$ is symmetric, and $S_m$ is symmetric CI 
if and only if $S_{m-1}$ is symmetric CI.
\end{lemma}
For our purpose the following Corollary of Lemma \ref{le1} is important.
\begin{corollary}\label{cor1}
Let a numerical semigroup $S_{m-1}\!=\!\la\delta_1,\ldots,\delta_{m-1}\ra$ be 
given and $a\in{\mathbb Z}_{>0}$, $a>1$, $d_m>m$, such that $\gcd(a,d_m)=1$, 
$d_m\in S_{m-1}$. Consider a semigroup $S_m\!=\!\!\la aS_{m-1},d_m\ra$. Then 
$S_m$ is symmetric (not CI) if and only if $S_{m-1}$ is symmetric (not CI).  
\end{corollary}
To utilize this construction we define the following property.
\begin{definition}\label{de1}
{\rm A symmetric (not CI) semigroup $S_m$ has {\em the property W} if there 
exists another symmetric (not CI) semigroup $S_{m-1}$ giving rise to $S_m$ 
by the construction, described in Corollary \ref{cor1}.}
\end{definition}
Note that in Definition \ref{de1} the semigroup $S_{m-1}$ does not necessarily
possess the property {\em W}. The symbol {\em W} stands for Kei-ichi Watanabe. 
The next Proposition distinguishes the minimal {\em edim} of symmetric CI and
symmetric (not CI) numerical semigroups with the property {\em W} .
\begin{proposition}\label{pr2}
A minimal {\em edim} of symmetric (not CI) semigroup $S_m$ with the property
{\em W} is $m=5$.
\end{proposition}
\begin{proof}
All numerical semigroups of $edim=2$ are symmetric CI, and all symmetric
numerical semigroups of $edim=3$ are CI \cite{herz70}. Therefore, according to 
Definition \ref{de1}, a minimal {\em edim} of symmetric CI semigroups with the 
property $W$ is $edim\!=\!3$. A minimal {\em edim} of numerical semigroups, 
which are symmetric (not CI), is $edim=4$. Therefore, according to Definition
\ref{de1} and Corollary \ref{cor1}, a minimal {\em edim} of symmetric (not CI)
semigroups with the property {\em W} is $edim=5$.
\end{proof}
According to Proposition \ref{pr2}, let us choose a symmetric (not CI) semigroup
of $edim=5$ with the property $W$ and denote it by ${\sf W}_5$ in order to 
distinguish it from other symmetric (not CI) semigroups ${\sf S}_5$ 
(irrespective of the $W$ property). Then the following containment holds: 
$\{{\sf W}_5\}\subset\{{\sf S}_5\}$. A minimal free resolution associated with 
semigroups ${\sf W}_5$ was described recently (\cite{bf12}, section 4), where 
the degrees of all syzygies were also derived (\cite{bf12}, Corollary 12), e.g.,
its 1st Betti number is $\beta_1=6$.
\begin{lemma}\label{le2}
Let two symmetric (not CI) semigroups ${\sf W}_5=\la a{\sf S}_4,d_5\ra$ and
${\sf S}_4=\la\delta_1,\delta_2,\delta_3,\delta_4\ra$  be given and $\gcd(a,d_5)
=1$, $d_5\in {\sf S}_4$. Let the lower bound $F_{5w}$ of the Frobenius number
$F({\sf W}_5)$ of the semigroup ${\sf W}_5$ be represented as, $F_{5w}=g_{5w}-
\left(a\sum_{j=1}^4\delta_j+d_5\right)$. Then
\bea
g_{5w}=a\left(\sqrt[3]{25\pi_4(\delta)}+d_5\right),\qquad\pi_4(\delta)=
\prod_{j=1}^4\delta_j.\label{h13}
\eea
\end{lemma}
\begin{proof}
Consider a symmetric (not CI) numerical semigroup ${\sf S}_4$ generated by four
integers (without the $W$ property), and apply the recent result \cite{fel015} 
on the lower bound $F_4$ of its Frobenius number, $F({\sf S}_4)$, to obtain
\bea
F({\sf S}_4)\ge F_4,\qquad F_4=h_4-\sum_{j=1}^4\delta_j,\qquad h_4=\sqrt[3]{25
\pi_4(\delta)}.\label{h14}
\eea
The following relationship between the Frobenius numbers $F({\sf W}_5)$ and 
$F({\sf S}_4)$ was derived in \cite{bra62}:
\bea
F({\sf W}_5)=aF({\sf S}_4)+(a-1)d_5.\label{h15}
\eea
Substituting two representations $F({\sf S}_4)\!=\!h-\sum_{j=1}^4\delta_j$ and 
$F({\sf W}_5)\!=\!g-a\sum_{j=1}^4\delta_j-d_5$ into equality (\ref{h15}), we  
obtain
\bea
g-a\sum_{j=1}^4\delta_j-d_5=ah-a\sum_{j=1}^4\delta_j+(a-1)d_5\qquad\rightarrow
\qquad g=a(h+d_5).\label{h16}
\eea
Comparing the last equality in (\ref{h16}) with the lower bound of $h$ in
(\ref{h14}), we arrive at (\ref{h13}).
\end{proof}
\begin{proposition}\label{pr3}
Let two symmetric (not CI) semigroups ${\sf W}_5=\la a{\sf S}_4,d_5\ra$ and
${\sf S}_4=\la\delta_1,\delta_2,\delta_3,\delta_4\ra$, where $\gcd(a,d_5)=1$, 
$d_5\in{\sf S}_4$, be given with the Hilbert series $H\left({\sf W}_5;z\right)$ 
according to (\ref{h1}). Then
\bea
g_{5w}>g_5.\label{h17}
\eea
\end{proposition}
\begin{proof}
Keeping in mind that $\beta_1=6$ (see \cite{bf12}, Corollary 12), we determine
$g_5$ for the semigroup ${\sf W}_5$ according to (\ref{h9}),
\bea
g_5=a\lambda(6)\sqrt[4]{\pi_4(\delta)d_5},\qquad\lambda(6)=4\sqrt[4]{5/8}\simeq
3.556,\nonumber
\eea
and calculate the following ratio:
\bea
\rho=\frac{g_{5w}}{g_5}=\frac1{\lambda(6)}\frac{\sqrt[3]{25\pi_4(\delta)}+d_5}
{\sqrt[4]{\pi_4(\delta)d_5}}=\frac1{\lambda(6)}\left(\sqrt[3]{25\eta}+\frac1
{\eta}\right),\qquad\eta=\sqrt[4]{\frac{\pi_4(\delta)}{d_5^3}}.\nonumber
\eea
The function $\rho(\eta)$ is positive for $\eta>0$ and has an absolute minimum
\bea
\rho(\eta_m)\simeq 1.1033,\qquad\eta_m\simeq 1.01943,\qquad\pi_4(\delta)\simeq
1.08 d_5^3.\nonumber
\eea
In other words, we arrive at $\rho(\eta_m)>1$, which proves the Proposition.
\end{proof}
We present ten symmetric (not CI) semigroups $A_j^6$ with the $W$ property
and the Betti number $\beta_1=6$. All semigroups $A_j^6$ are built according to 
Lemma \ref{le1} and based on symmetric (not CI) semigroups ${\sf S}_4$, 
generated by four integers \cite{fel017}. We give a comparative Table 2 for 
their $g$, lower bounds $g_{5w}$ and $g_5$, and the parameter $\eta$.
\bea
&&A_1^6=\!\la 14,15,16,18,26\ra,\quad A_2^6=\!\la 20,21,24,27,39\ra,\quad
A_3^6=\!\la 10,11,12,14,16\ra,\nonumber\\
&&A_4^6=\!\la 14,16,17,18,26\ra,\quad A_5^6=\la 16,21,26,30,34\ra,\quad
A_6^6=\la 23,24,39,45,51\ra,\nonumber\\
&&A_7^6=\la 35,40,41,45,65\ra,\quad A_8^6=\la 302,305,308,314,316\ra,
\nonumber\\
&&A_9^6=\la 302,308,314,315,316\ra,\quad A_{10}^6=\la 453,462,469,471,474\ra,
\nonumber
\eea
$$
\begin{array}{|c||c|c|c|c|c|c|c|c|c|c|c|c|}\hline
{\sf W}_5 & A_1^6 & A_2^6 & A_3^6 & A_4^6 & A_5^6 & A_6^6 & A_7^6 & A_8^6 
& A_9^6 & A_{10}^6\\ \hline\hline
g & 142 & 228 & 92 & 146 & 218 & 333 & 485 & 9120 & 9140 & 14172 \\ \hline
g_{5w}& 139.4 & 224.1 & 91.5 & 143.4 & 216.4 & 330.6 & 478.6 & 5478.1 & 5498.1
& 8709.2\\ \hline
g_5& 125.9 & 203.0 & 82.9 & 129.9 & 194.3 & 298.2 & 404.8 & 4606.8 & 4644.1 &
7695.0\\ \hline
\eta& 1.180 & 0.951 & 1.060 & 1.075 & 1.301 & 1.215 & 0.555 & 2.123 & 2.073 &
1.538 \\ \hline
\end{array}
$$

\begin{center}
Table 2. The largest degree $g$ of syzygies, its lower bounds $g_5$ and $g_{5w}$
and the parameter $\eta$ for symmetric (not CI) semigroups ${\sf W}_5$ with the 
Betti number $\beta_1=6$.
\end{center}
\section*{Acknowledgement}
The author is thankful to M. Delgado for providing him with various examples of 
symmetric semigroups generated by five elements which were found numerically 
with the help of the package "NumericalSgps" \cite{dgm00}. The author is 
grateful to the managing editor for valuable suggestions to improve the
manuscript.

\end{document}